\def\L{{\mathcal L}}
\def\C{{\mathcal C}}
\def\A{{\mathcal A}}
\def\N{\mathbb N}
\def\Z{\mathbb Z}
\def\Q{\mathbb Q}
\def\pf{\begin{proof}}
\def\pfk{\end{proof}}
\documentclass[a4paper,11pt,reqno,oneside]{article}

\usepackage{amsmath}
\usepackage{amsthm}
\usepackage{amssymb}

\usepackage{enumerate}
\usepackage{epic}

\newtheorem{thm}{Theorem}[section]
\newtheorem{prop}[thm]{Proposition}
\newtheorem{lem}[thm]{Lemma}
\newtheorem{coro}[thm]{Corollary}
\begin{document}

\begin{center}

{\LARGE\bf Relation between powers of factors and recurrence\\[4mm]
function characterizing Sturmian words}

\vskip0.8cm

{\large Z. Mas\'akov\'a\footnote{corresponding author} \quad and
\quad E. Pelantov\'a}

\vskip0.5cm

{Doppler Institute \& Department of Mathematics\\
FNSPE, Czech Technical University\\
Trojanova 13, 120 00 Praha 2, Czech Republic\\

e-mails: zuzana.masakova@fjfi.cvut.cz,
edita.pelantova@fjfi.cvut.cz}

\vskip 0.5cm

\begin{abstract}
In this paper we use the relation of the index of an infinite
aperiodic word and its recurrence function to give another
characterization of Sturmian words. As a byproduct, we give a new
proof of theorem describing the index of a Sturmian word in terms
of the continued fraction expansion of its slope. This theorem was
independently proved in~\cite{Carpi} and~\cite{DaLenzNecelo}.
\end{abstract}

\end{center}

\section{Introduction}

Sturmian words constitute the most studied example of aperiodic
infinite words. For the first time they appeared in the paper of
Morse and Hedlund in 1938~\cite{morse}. But even after 70 years of
extensive research, Sturmian words continue to attract attention
of numerous mathematicians and newly also computer scientists. The
appeal of Sturmian words stems in that they appear in various
contexts. This is also why Sturmian words are often hidden under
different titles: cutting sequences, Beatty sequences, mechanical
words, etc. The beauty of Sturmian words consists in the abundance
of equivalent definitions. Already Morse and Hedlund
in~\cite{hedlund} show that Sturmian words can be characterized by
the so-called balance property. The reference~\cite{lothaire}
contains a nice exposition on diverse definitions of Sturmian
words. The most recent ones, which~\cite{lothaire} does not
mention, are characterization of Sturmian words using return words
given by Vuillon~\cite{vuillon} (for less technical proof
see~\cite{wolfi}), characterization using the number of
palindromes of given length given in~\cite{Droubay} and yet
another characterization by Richomme~\cite{richome}.

The aim of this paper is to give another equivalent definition of
Sturmian words. Our characteristics puts into relation the
recurrence function and the index of an infinite word $u$.
Flagrant similarity between formulas for recurrence quotient and
index of a Sturmian word was noted already
in~\cite{mirrorformula,initialpowers,Carpi}.

The recurrence function $R$ associates to every $n\in\N$ the
minimal length $R(n)\in\N$ such that arbitrary segment of the
infinite word $u$ of length $R(n)$ contains all factors of $u$ of
length $n$. This function has been studied already by Hendlund and
Morse, who gave an explicit formula for $R(n)$ for an arbitrary
Sturmian word $u$ and determined the so-called recurrence
quotient, $\limsup_{n\to\infty}{R(n)}/{n}$.
 On the other hand, the index of an infinite word $u$
describes the maximal repetition of a factor of $u$. The study of
the index of infinite words is considerably younger, nevertheless,
in the last decade very intense, especially due to applications in
spectral theory for corresponding Schr\"odinger
operators~\cite{damanik}.

Repetitions in the most prominent Sturmian word, namely the
Fibonacci word, were studied in~\cite{karhumaki}. More general
results about index of Sturmian words can be found
in~\cite{berstel2,initialpowers,CaoWen,JuPiPowers,MignosiPi,Mignosi,Vandeth}.
The complete solution to the problem was given independently by
Carpi and de Luca in~\cite{Carpi} and by Damanik and Lenz
in~\cite{DaLenzNecelo}.

The paper is organized as follows. In Section~\ref{sec:preli} we
introduce all necessary notions. Section~\ref{sec:3} contains the
proof of the main result of the paper, namely the following
theorem.

\begin{thm}\label{t:A}
A uniformly recurrent infinite word $u$ is Sturmian if and only if
there exist infinitely many factors $w$ of $u$ such that
$$
R(|w|) = |w|\,{\rm ind}(w)+1\,.
$$
\end{thm}

Notation $|w|$ stands for the length of the factor $w$, and ${\rm
ind}(w)$ is the maximal rational exponent $r$ such that $w^r$ is a
factor of $u$.

It was pointed to us that already from~\cite{Carpi} one can
extract that Sturmian words satisfy the above equality for
infinitely many factors. Their proof uses the explicit formula for
recurrence function from~\cite{hedlund}. Our proof relies on
Vuillon's description of Sturmian words by return words and avoids
manipulation with continued fraction of the slope of the Sturmian
word. Our theorem moreover states that Sturmian words are the only
having the above property.

With the help of Theorem~\ref{t:A}, one can derive the upper bound
on the index of $u$ (Section~\ref{sec:4}). In Section~\ref{sec:5}
we prove that the bound is in fact reached. For the construction
of factors of $u$ with large repetition we use the knowledge of
Sturmian morphisms, i.e.\ morphisms preserving the family of
Sturmian words, as described in~\cite{seebold}.
Sections~\ref{sec:4} and~\ref{sec:5} thus represent an alternative
proof of the result of~\cite{Carpi} and~\cite{DaLenzNecelo}.

\section{Preliminaries}\label{sec:preli}

An alphabet $\mathcal A$ is a finite set of symbols, called
A word $w$ of length $|w|=n$ is a concatenation of $n$ letters.
The number of letters $X$ occurring in the word $w$ is denoted by
$|w|_X$. ${\mathcal A}^*$ is the set of all finite words over the
alphabet $\A$ including the empty word $\epsilon$. Equipped with
the operation of concatenation, it is a monoid. We define also
infinite words $u=(u_n)_{n\in\N}\in\A^{\N}$.

A finite word $v\in{\A}^*$ is called a {\em factor} of a word $w$
(finite or infinite), if there exist words $w^{(1)},w^{(2)}$ such
that $w=w^{(1)}vw^{(2)}$. If $w^{(1)}=\epsilon$, then $v$ is said
to be a prefix of $w$, if $w^{(2)}=\epsilon$, then $v$ is a suffix
of $w$. The set of all factors of length $n$ of an infinite word
$u$ is denoted by $\L_n(u)$, the set of all factors of $u$ is
called the {\em language} of $u$ and denoted by $\L(u)$.

The mapping $\C: n\mapsto \#\L_n(u)$ is called the {\em
complexity} of the infinite word $u$. For determining the
complexity of an infinite word one uses the so-called special
factors. A factor $w\in\L(u)$ is called {\em left special}, if
there exist letters $A,B\in\A$, $A\neq B$, such that both $Aw$ and
$Bw$ belong to $\L(u)$. Similarly, one defines {\em right special}
factors. A factor of $u$ is called bispecial, if it is in the same
time right special and left special. Every eventually periodic
word has bounded complexity. For aperiodic words, one has for all
$n\in\N$ that $\C(n)\geq n+1$. Infinite words, for which equality
holds for all $n\in\N$, i.e.\ aperiodic words with minimal
complexity, are called {\em Sturmian words}. Directly from the
definition one can derive that in the language of a Sturmian word
$u$ one has exactly one left special and exactly one right special
factor of each length, and Sturmian words are characterized by
this property.

Sturmian words are obviously defined over a binary alphabet, say
$\{A,B\}$. The densities of letters $A$, $B$ in a Sturmian word
$u=(u_i)_{i\in\N}$ are well defined,
$$
\varrho(A)=\lim_{n\to\infty}\frac{|u_0\cdots
u_{n-1}|_A}{n}=\alpha\,,\qquad
\varrho(B)=\lim_{n\to\infty}\frac{|u_0\cdots
u_{n-1}|_B}{n}=1-\alpha\,,
$$
for some $\alpha\in(0,1)$. In fact, the language of a Sturmian
word $u$ depends only on the parameter $\alpha$, which is also
called the {\em slope} of $u$. For a given $\alpha$, one can
construct all Sturmian words with the slope $\alpha$ for example
as codings of different orbits under an exchange of two intervals.
Let $\alpha\in(0,1)$ be an irrational number. Denote $I=[0,1)$
(resp. $I=(0,1]$) and $I_A=[0,\alpha)$, $I_B=[\alpha,1)$ (resp.
$I_A=(0,\alpha]$, $I_B=(\alpha,1]$). The mapping $T:I\mapsto I$
given by the prescription
$$
T(x)=\left\{\begin{array}{ll} x+1-\alpha & \hbox{ for } x\in
I_A\,,\\
x-\alpha & \hbox{ for } x\in I_B\,,
\end{array}\right.
$$
is called an exchange of two intervals with {slope} $\alpha$. For
an arbitrary $x_0\in I$ we define an infinite word
$(u_n)_{n\in\N}$ by
\begin{equation}\label{eq:7}
u_n=X\in\{A,B\} \quad\hbox{if}\quad T^n(x_0)\in I_X\,.
\end{equation}
It is known that the set of Sturmian words coincides with the set
of infinite words given by the prescription~\eqref{eq:7}. Since we
assume that the slope is irrational, the language of a Sturmian
word does not depend on the choice of the initial point $x_0$, but
only on $\alpha$. Due to the symmetry $\alpha \leftrightarrow
1-\alpha$, studying the language of a Sturmian word, one can
consider without loss of generality only parameters
$\alpha>\frac12$. From the exchange of intervals is not difficult
to see that with such an assumption, $\varrho(A)>\varrho(B)$ and,
in fact, the Sturmian word can be viewed as composed by blocks of
the form $A^k$, $A^{k+1}$, with
$k=\lfloor\frac{\alpha}{1-\alpha}\rfloor$, separated by single
letters $B$.

In this paper we study repetition of factors in Sturmian words. We
say that a word $v$ is a {\em power} of a word $w$, if
$|v|\geq|w|$ and $v$ is a prefix of the periodic word $www\cdots$.
We write $v=w^r$ where $r=|v|/|w|$. The index of a word $w$ in an
infinite word $u$ is defined by
\begin{equation}\label{e:indfactor}
{\rm ind}(w)= \sup\{r\in\Q\mid w^r\in\L(u)\}\,.
\end{equation}
A power $v$ of $w$ with maximal $r$ is called a {\em maximal
repetition} of $w$. We have thus $v=w^{{\rm ind}(w)}$. From what
it was said above, it is clear that in a Sturmian word with slope
$\alpha>\frac12$, one has
\begin{equation}\label{e:indexpismene}
{\rm ind}(B)=1 \qquad\hbox{ and }\qquad {\rm
ind}(A)=\Bigl\lfloor\frac{\alpha}{1-\alpha}\Bigr\rfloor+1\,.
\end{equation}

Taking supremum of indices over all factors of an infinite word
$u$, one obtains an important characteristics of $u$, the
so-called {\em index} of $u$. Formally,
\begin{equation}\label{e:indword}
{\rm ind}(u)= \sup\{{\rm ind}(w)\mid w\in\L(u)\}\,.
\end{equation}

It turns out that for the study of index of Sturmian words, the
notion of return words and recurrence function is important. A
{\em return word} of a factor $w$ of an infinite word $u$ is a
factor $v\in\L(u)$ such that $vw\in\L(u)$, $w$ is a prefix of $vw$
and the factor $w$ occurs in $vw$ exactly twice. The factor $vw$
is often called a {\em complete return word} of $w$. The set of
return words of a factor $w$ is denoted by ${\rm Ret}(w)$. If the
set ${\rm Ret}(w)$ is finite for any factor $w$ of an infinite
word $u$, then $u$ is said to be uniformly recurrent. In fact, it
means that distances between consecutive occurrences of a given
factor are bounded. Let us mention that for a uniformly recurrent
word $u$ the supremum in~\eqref{e:indfactor} is always reached, as
will be explained later, and therefore the notion of index of $u$
in~\eqref{e:indword} has sense. For a uniformly recurrent infinite
word $u$ we define a mapping $R:\N\mapsto\N$ by the prescription
\begin{equation}\label{eq:recurrencefunction}
R(n):= -1+\max\bigl\{|vw|\,\bigm|\, v\in{\rm Ret}(w), \,
w\in\L_n(u)\bigr\}\,,
\end{equation}
i.e.\ $R(n)+1$ is equal to the maximum of lengths of a complete
return word over all factors of length $n$. It is not difficult to
see that an arbitrary segment of the infinite word $u$ of length
$R(n)$ contains all factors of the word $u$ of length $n$. Formally,
we have
\begin{equation}\label{eq:vsechnyfaktory}
\L_n(u)=\{u_iu_{i+1}\cdots u_{i+n-1}\mid k\leq i\leq
k+R(n)-n+1\}\,,\quad \hbox{for all }\  k\in\N\,.
\end{equation}
Moreover, the number $R(n)$ is the smallest possible, so
that~\eqref{eq:vsechnyfaktory} remains valid. The mapping $R(n)$
is called the {\em recurrence function} of the infinite word $u$.


\section{Recurrence function and index}\label{sec:3}

Our aim is to find relation between the recurrence function (well
defined for uniformly recurrent words) and the index of aperiodic
words. We first show that index of every factor in an aperiodic
uniformly recurrent word is finite, and we then determine a lower
bound on the recurrence function.

\begin{prop}\label{p:2}
Let $u$ be an aperiodic uniformly recurrent word. Then for every
factor $w\in\L(u)$ we have $\ {\rm ind}(w)<+\infty\ $ and
\begin{equation}\label{e:hlavnivztah}
R(|w|)\ \geq \ |w|\, {\rm ind}(w) + \C(|w|) - |w|\,.
\end{equation}
\end{prop}

\pf Let $w=w_1\cdots w_n$ be a factor of $u$. We first show that
${\rm ind}(w)$ is finite. Without loss of generality, let ${\rm
ind}(w)\geq 2$. Obviously, all factors of the form $w_i\cdots
w_nw_1\cdots w_{i-1}$ for any $1\leq i\leq n$ belong to $\L(u)$.
(Such factors are called conjugates of $w$.) Since $\C(n)\geq
n+1$, there exists a factor $w'$ which is not conjugate of $w$. If
$\L(u)$ contained factors $w^k$ for all $k\in\N$, then distances
between consecutive occurrences of $w'$ would be unbounded, which
would contradict uniform recurrence of $u$. Therefore ${\rm
ind}(w)<+\infty$.

Let now $v$ be a maximal repetition of $w$. We prolong $v$ to a
factor $vv'\in\L(u)$ so that $vv'$ contains all $\C(|w|)$ factors
of $u$ of length $|w|$, but none of prefixes of $vv'$ satisfies
it. Since $v$ has at most $|w|$ factors of length $|w|$, (namely
the conjugates of $w$), we must have $|v'|\geq \C(|w|)-|w|$. From
the definition of the recurrence function, we have
$$
R(|w|)\geq |vv'|\geq |v| + \C(|w|)-|w|\,.
$$
As $v=|w|\,{\rm ind}(w)$, the proof is complete.
 \pfk

Note that in particular, for a Sturmian word $u$ one has $R(|w|)\
\geq \ |w|\, {\rm ind}(w) + 1$ for every factor $w$ of $u$. The
following proposition states, that if equality is reached for
infinitely many factors $w$ of an aperiodic word $u$, then $u$ is
Sturmian.

\begin{prop}\label{p:1}
Let $u$ be an aperiodic uniformly recurrent infinite word. If
there exist infinitely many factors $w\in\L(u)$ such that $R(|w|)\
= \ |w|\, {\rm ind}(w) + 1$, then $u$ is a Sturmian word.
\end{prop}

\pf Using the assumption of the proposition
and~\eqref{e:hlavnivztah}, there exist infinitely many factors $w$
of $u$ such that $\C(|w|)\leq |w|+1$, i.e.\ for infinitely many
$n\in\N$ we have $\C(n)\leq n+1$. The complexity of an aperiodic
word is a strictly increasing function and $\C(1)\geq 2$. This
implies that $\C(n)= n+1$ for all $n$ and $u$ is therefore
Sturmian. \pfk

In order to show the opposite implication to that of
Proposition~\ref{p:1}, we need to cite a nice result of
Vuillon~\cite{vuillon} which characterizes Sturmian words using
return words. He shows that a binary infinite word $u$ is Sturmian
if and only if every factor of $u$ has exactly two return words.
For every factor $w$ of a Sturmian word $u$ thus exist two finite
words $r_0(w)$, $r_1(w)$ such that the suffix of $u$ starting with
the first occurrence of $w$ can be written as an infinite
concatenation of blocks $r_0(w)$ and $r_1(w)$, i.e.
$$
u = p\,r_{i_0}(w)r_{i_1}(w)r_{i_2}(w)r_{i_3}(w)\cdots\,,
$$
where $p$ is a prefix of $u$ and $i_0,i_1,i_2,i_3,\cdots
\in\{0,1\}$. We can therefore define the so-called derivated word
$v=(v_n)_{n\in\N}$ over the alphabet $\{0,1\}$ by the prescription
$v_n=i_n$, coding the order of the blocks $r_0(w)$, $r_1(w)$ in
the infinite concatenation. We could now study return words of
factors of the newly defined infinite word $v$. However, since
return words of factors of the derivated word are in one-to-one
correspondence with return words of factors in the original
infinite word (see~\cite{durand}), we deduce that every factor of
$v$ has again exactly two return words, and thus is itself
Sturmian.

It is obvious that for finding factors $w$ with the maximal index
in the infinite word, we can limit our consideration to primitive
factors $w$, i.e.\ such that $w\neq z^k$ for any $z\in\L(u)$ and
any $k\in\N$, $k\geq 2$.

\begin{prop}\label{t:rekindex}
Let $u$ be a Sturmian word and let $w\in\L(u)$ be a primitive
factor such that $ww\in\L(u)$, and, moreover, let it have the
maximal index among all factors of $u$ of length $n$ with the
above properties. Then
$$
R(n)=n\,{\rm ind}(w)+1\,.
$$
\end{prop}

\pf Let $k=[{\rm ind}(w)]$ and $\theta=\{{\rm ind}(w)\}$. Then $w$
can be written as $w=w_1w_2$ where $|w_1|=\theta n$ and the
maximal repetition of $w$ is the word
$$
\underbrace{(w_1w_2)(w_1w_2)\cdots(w_1w_2)}_{k\ \hbox{\scriptsize
times}}w_1 \in\L(u)\,.
$$
Let us find $X,Y\in\{A,B\}$ such that
\begin{equation}\label{eq:bispecial}
Xw_1w_2\cdots w_1w_2w_1Y\in\L(u)\,.
\end{equation}
Since ${\rm ind}(w)=k+\theta$ is the greatest power such that
$w^{k+\theta}\in\L(u)$, the letter $Y$ is not a prefix of $w_2$.
Since $w$ is a primitive word with the greatest index in
$\L_n(u)$, the letter $X$ is not a suffix of $w_2$. This, together
with the fact that $k\geq 2$, means that $w_1w_2=w$ is a left
special factor and $w_2w_1=:w'$ is a right special factor. A
Sturmian word has exactly one left special and one right special
factor of each length.

Let us consider the Rauzy graph $\Gamma_n$ of $u$. The set of
vertices of $\Gamma_n$ is equal to $\L_n(u)$ and the set of its
edges to $\L_{n+1}(u)$. The Rauzy graph $\Gamma_n$ of a Sturmian
word thus has ${n+1}$ vertices and $n+2$ edges. An edge
$e\in\L_{n+1}(u)$ starts in a vertex $v\in\L_{n}(u)$ and ends in
$v'\in\L_{n}(u)$ if $v$ is a prefix and $v'$ a suffix of $e$. An
arbitrary factor $u$ of length $m\geq n$ in the language of the
infinite word $u$ can be viewed as a path of length $m-n$ in the
graph $\Gamma_n$ starting in the vertex corresponding to the
prefix and ending in the vertex corresponding to the suffix of $u$
of length $n$.

Since $w\in\L_n(u)$, $ww\in\L(u)$ and $w$ is primitive, there
exists a cycle $C$ of length $n$ in the graph $\Gamma_n$
containing the factor $w$. Let us denote the vertices of the cycle
$C$ by $v^{(0)}=w$, $v^{(1)}$, \dots, $v^{(n-1)}$. Since
$\Gamma_n$ has $n+1$ vertices, only one of them is missing in $C$.
Let us denote it by $v^{(n)}$. Recall that $w$ is the only left
special factor in $\L_n(u)$, and thus the only vertex in
$\Gamma_n$ with indegree 2. Similarly, $w'$ is the only right
special factor in $\L_n(u)$ and thus the only vertex in $\Gamma_n$
with outdegree 2. Since $\Gamma_n$ is a strongly connected graph,
an edge must go from the vertex $v^{(n)}$ to the cycle $C$ and an
edge from the cycle $C$ to the vertex $v^{(n)}$. Thus $w'=v^{(s)}$
for some $0\leq s\leq n-1$. Relation~\eqref{eq:bispecial} implies
that the edge from $v^{(s)}$ to $v^{(n)}$ is $w_2w_1Y$ and the
edge from $v^{(n)}$ to $v^{(0)}$ is $Xw_1w_2$. The Rauzy graph
$\Gamma_n$ is thus of the following form.

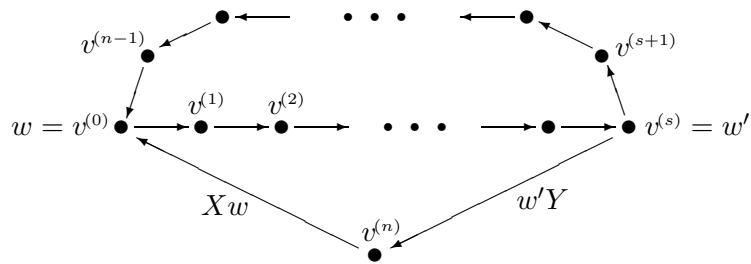
\begin{center}
\setlength{\unitlength}{1pt}
\begin{picture}(362,109)
\put(44,50){$w=v^{\scriptscriptstyle\!(0)}$}
\put(85,53){\circle*{5}} \put(115,53){\circle*{5}}
\put(110,58){$v^{\scriptscriptstyle\!(1)}$}
\put(145,53){\circle*{5}}
\put(140,58){$v^{\scriptscriptstyle\!(2)}$}
\put(185,53){\circle*{3}} \put(195,53){\circle*{3}}
\put(205,53){\circle*{3}} \put(245,53){\circle*{5}}
\put(275,53){\circle*{5}} \put(90,53){\vector(1,0){20}}
\put(120,53){\vector(1,0){20}} \put(150,53){\vector(1,0){20}}
\put(220,53){\vector(1,0){20}} \put(250,53){\vector(1,0){20}}
\put(281,50){$v^{\scriptscriptstyle\!(s)}=w'$}
\put(175,10){$v^{\scriptscriptstyle\!(n)}$}
\put(180,5){\circle*{5}} \put(269,49){\vector(-2,-1){83}}
\put(174,7.5){\vector(-2,1){83}}
\put(265,80){\circle*{5}}
\put(270,80){$v^{\scriptscriptstyle\!(s+1)}$}
\put(273.5,57.5){\vector(-1,3){6.3}}
\put(237,94.6){\circle*{5}}
\put(261.5,83){\vector(-2,1){20}}
\put(232,94.6){\vector(-1,0){20}}
\put(95,80){\circle*{5}}
\put(70,81){$v^{\scriptscriptstyle\!(n-1)}$}
\put(93.5,75.5){\vector(-1,-3){6.3}}
\put(123,94.6){\circle*{5}}
\put(118.5,92.6){\vector(-2,-1){19}}
\put(148,94.6){\vector(-1,0){20}}
\put(170,94.6){\circle*{3}}
\put(180,94.6){\circle*{3}}
\put(190,94.6){\circle*{3}}
\put(115,20){$Xw$} \put(233,22){$w'Y$}
\end{picture}
\end{center}

Let us consider the return words of $w$. Since $ww\in\L(u)$, one
of the return words of $w$ is $r_0(w)=w$, the complete return word
is $ww$ and the corresponding path in the Rauzy graph is the cycle
$C$. We denote the other return word of $w$ by $r_1(w)$. From the
structure of the graph $\Gamma_n$ it follows that the complete
return word $r_1(w)w$ corresponds to the cycle $C'$ given by
vertices $v^{(0)}$, $v^{(1)}$, \dots, $v^{(s)}$, $v^{(n)}$.

As we have already mentioned, the order of the blocks $r_0(w)$,
$r_1(w)$ is given by the derivated word over the alphabet
$\{0,1\}$, which is Sturmian. Since $(r_0(w))^k=w^k\in\L(u)$, for
$k=\lfloor {\rm ind}(w)\rfloor\geq 2$, the derivated word has
blocks $0^k$, $0^{k-1}$ separated by single letters 1. As a
consequence, among all factors of length $n$, it is $v^{(n)}$
which has the longest complete return word, namely of the form
$$
X\underbrace{ww\cdots w}_{k\ \hbox{\scriptsize
times}}w_1Y\,.
$$
From the definition~\eqref{eq:recurrencefunction} it follows that
$$
R(n) = -1+|w^{k+\theta}|+2 = 1+(k+\theta)n\,,
$$
which completes the proof.
\pfk

\pf[Proof of Theorem~\ref{t:A}] In order to comlpete the proof of
Theorem~\ref{t:A}, we have to show that there exist infinitely
many primitive factors $w$ with index at least 2. For the
construction of such factors we make use of bispecial factors. Let
$b$ be a bispecial factor in $\L(u)$. Denote by $n$ its length,
$n:=|b|$ and by $r_0(b)$, $r_1(b)$ its return words. From the
Rauzy graph $\Gamma_n$ it follows that the two return words of $b$
are given by the two cycles in $\Gamma_n$, which have $b$ as the
only common vertex. Therefore $|r_0(b)|+|r_1(b)|=n+2$. Without
loss of generality, let $b$ contain both letters. Then
$|r_i(b)|\geq 2$. At least for one of the return words, say
$r_0(b)$, it holds that $n/2<|r_0(b)|\leq n$, and therefore
$r_0(b)$ is a prefix of $b$. It follows that the complete return
word $r_0(b)b\in\L(u)$ has as its prefix $r_0(b)r_0(b)$. Moreover,
a return word of an arbitrary factor of any uniformly recurrent
word is primitive. Thus we can take $r_0(b)$ for the desired
factor $w$. Since there are infinitely many bispecial factors $b$,
we can construct infinitely many primitive factors with index
$\geq 2$ and length $\geq \frac{|b|}{2}$.
 \pfk

\section{Upper bound on index of Sturmian words}\label{sec:4}

In this section we mention the consequences of
Proposition~\ref{t:rekindex}, which puts into relation the
recurrence function and index of factors of a Sturmian word. In
particular, we can very easily derive the upper bound on the index
of a Sturmian word, which constitutes an alternative proof for the
result of Damanik and Lenz~\cite{DaLenzNecelo}. The bound depends
on the continued fraction expansion of the slope of the Sturmian
word.

Recall the notion of continued fraction. To every irrational
$\beta\in(0,1)$ one associates the continued fraction
$\beta=[0,b_1,b_2,\dots]$, where $b_i\in\Z$, $b_i\geq 1$.
Obviously, if $\beta>\frac12$, then $b_1=1$. The convergents of
$\beta$ form a sequence of fractions $(\frac{p_n}{q_n})$,
$$
\frac{p_1}{q_1}=\frac1{b_1}\,, \qquad
\frac{p_2}{q_2}=\cfrac1{b_1+\cfrac1{b_2}}\,, \qquad
\frac{p_3}{q_3}=\cfrac1{b_1+\cfrac1{b_2+\cfrac1{b_3}}}\,,\qquad\dots
$$
We have $p_n$ coprime to $q_n$ and $\lim_{n\to\infty}\frac{p_n}{q_n}
= \beta$.

It is known that the denominators $q_n$ of convergents of $\beta$
satisfy the recurrence $$
q_N=b_Nq_{N-1}+q_{N-2}
$$
with initial values $q_{-1}=0$, $q_0=1$. Denoting the matrix
$M_c:=\bigl(\begin{smallmatrix}c&1\\1&0\end{smallmatrix}\bigr)$,
then the recurrence can be rewritten as
$$
(q_N,q_{N-1}) = (q_{N-1},q_{N-2}) M_{b_N}\,,
$$
and by repetition, we obtain
$$
(q_N,q_{N-1}) = (1,0) M_{b_1} M_{b_2} \cdots M_{b_N}
$$
In order to extract the component $q_N$, it suffices to multiply the
latter from the right by the vector $\binom{1}{0}$. We obtain
\begin{equation}\label{eq:01}
q_N=(1,0) M_{b_1} M_{b_2} \cdots M_{b_N}\textstyle{\binom{1}{0}} =
(1,0) M_{b_N} \cdots M_{b_2} M_{b_1}\textstyle{\binom{1}{0}}\,,
\end{equation}
where we have used that equality must hold also for the transpose
$q_N^T=q_N$ and $M_c^T=M_c$ for all $c\in\N$.

For the derivation of the lower bound on the index of Sturmian
words we use an old result on recurrence function of Sturmian
words given in~\cite{hedlund}.

\begin{thm}[\cite{hedlund}]\label{t:hedlund}
Let $u$ be a Sturmian word with slope $\alpha$. Denote by
$q_0,q_1,q_2,\dots$ the denominators of the convergents of $\alpha$.
Then for every $n\in\N$,
$$
R(n)=q_{N+1}+q_N+n-1\,,\qquad\hbox{where $N$ is such that }\
q_N\leq n < q_{N+1}\,.
$$
\end{thm}

Substituting into Proposition~\ref{t:rekindex}, one obtains an
easy proof of the following result. Similar derivation one can
find in~\cite{Carpi}.

\begin{coro}\label{t:upperbound}
Index of every factor of a Sturmian word $u$ with the slope
$\alpha=[0,1,a_2,a_3,\dots]$ is bounded by
$$
\sup\Big\{\,2+a_{N+1}+\frac{q_{N-1}-2}{q_N}\;\Big|\; N\geq
1\,\Big\}\,,
$$
where $q_N$ are the denominators of the convergents of $\alpha$.
\end{coro}

\pf Obviously, it suffices to consider only factors $w$ satisfying
assumptions of Theorem~\ref{t:rekindex}. Let $|w|=n$ and let
$q_N\leq n<q_{N+1}$. Using Proposition~\ref{t:rekindex} and
Theorem~\ref{t:hedlund}, we have
$$
n\,{\rm ind}(w)+1\ =\ R(n)\ =\ q_{N+1} + q_N +n-1\,.
$$
Therefore
$$
q_N\bigl({\rm ind}(w)-1\bigr)\ \leq\ n\bigl({\rm ind}(w)-1\bigr)\
=\ q_{N+1} + q_N -2\ =\ (a_{N+1}+1)q_N + q_{N-1} -2\,,
$$
and consequently
$$
{\rm ind}(w) \ \leq \ 2+a_{N+1}+\frac{q_{N-1}-2}{q_N}\,.
$$
\pfk

\section{Sturmian morphisms and factors with maximal
index}\label{sec:5}

In this section we provide a lower bound on the index of a
Sturmian word $u$ of slope $\alpha$. Obviously, ${\rm ind}(u)\geq
a_2+1$, since $\lfloor\frac{\alpha}{1-\alpha}\rfloor$ in the
formula~\eqref{e:indexpismene} for the index of the letter $A$ is
equal to the coefficient $a_2$ of the continued fraction of
$\alpha$. The idea for construction of factors with large index in
a Sturmian word $u$ stems in application of specific Sturmian
morphisms. Since application of a morphism preserves repetitions,
it suffices to know how the chosen morphism changes the slope of
the Sturmian word. Let us recall the necessary facts.

A morphism over the alphabet $\{A,B\}$ is a mapping
$\varphi:\{A,B\}^*\mapsto \{A,B\}^*$ satisfying
$\varphi(w_1w_2)=\varphi(w_1)\varphi(w_2)$. Obviously, a morphism
is uniquely determined by $\varphi(A)$, $\varphi(B)$. The
incidence matrix of a morphism $\varphi$ is given by
$$
M_\varphi= \Biggl(\!\begin{array}{cc}|\varphi(A)|_A&|\varphi(A)|_B\\
|\varphi(B)|_A&|\varphi(B)|_B\end{array}\!\Biggr)
$$
The action of a morphism can be naturally extended to infinite words by
$$
\varphi(u_0u_1u_2\cdots) = \varphi(u_0) \varphi(u_1) \varphi(u_2)\cdots
$$
It is easy to show that for the number of letters in the image of a word $w$, one has
\begin{equation}\label{eq:p3}
\bigl(|\varphi(w)|_A,|\varphi(w)|_B\bigr) = \bigl(|w|_A,|w|_B\bigr) M_\varphi\,.
\end{equation}
From that, we can deduce the following fact for the densities of letters in an infinite word $u$.
If $\varrho(A)$, $\varrho(B)$ are the densities in $u$, than the densities in the word
$u'=\varphi(u)$ are $\varrho'(A)$, $\varrho'(B)$, where
\begin{equation}\label{eq:p2}
\bigl(\varrho'(A),\varrho'(B)\bigr) = {\it const.}\,\bigl(\varrho(A),\varrho(B)\bigr) M_\varphi\,,
\end{equation}
and ${\it const.}$ is chosen so that $\varrho'(A)+\varrho'(B)=1$.

A morphism $\varphi$ is called Sturmian, if $\varphi(u)$ is a Sturmian word for every Sturmian word $u$.
Obviously, the set of Sturmian morphisms equipped with the operation of composition is a monoid, denoted
by ${\it St}$. It is known~\cite{seebold} that the monoid ${\it St}$ has three generators, namely
\begin{equation}\label{eq:p1}
\psi_1: \begin{array}{rcl}
A&\mapsto& AB\\
B&\mapsto& B
\end{array}
\qquad
\psi_2: \begin{array}{rcl}
A&\mapsto& BA\\
B&\mapsto& B
\end{array}
\qquad
E: \begin{array}{rcl}
A&\mapsto& B\\
B&\mapsto& A
\end{array}
\end{equation}

Consider a Sturmian word with slope $\beta\in(\frac12,1)$ whose continued fraction is of the form
$\beta=[0,1,b_2,b_3,\dots]$.
For $c\in\N$, we shall study the action of the morphism
\begin{equation}\label{eq:101}
\varphi : \begin{array}{rcl}
A&\mapsto& A^cB\\
B&\mapsto& A
\end{array}
\end{equation}
on the Sturmian word $u$ with slope $\beta$. The morphism
$\varphi$ is a Sturmian morphism; it is a composition of the
generators~\eqref{eq:p1} of the Sturmian monoid, namely
$\varphi=E\psi_2^c$. The corresponding
incidence matrix is $M_{\varphi}=M_c=\bigl(\begin{smallmatrix}c&1\\
1&0\end{smallmatrix}\bigr)$, as defined in the Preliminaries.
Consequently, the infinite word $\varphi(u)$ is also Sturmian,
i.e.\ there exists an irrational $\beta'$ such that
$u':=\varphi(u)$ is a Sturmian word with slope $\beta'$. According
to~\eqref{eq:p2}, the densities of letters $a,b$ in the word $u'$
satisfy
$$
(\beta',1-\beta')\ =\ {\it const.}\ (\beta,1-\beta)\ \Bigl(\!\begin{array}{cc}c&1\\
1&0\end{array}\!\Bigr)\,.
$$
Therefore $\beta'=\frac{c\beta+1-\beta}{c\beta+1}$. It is not difficult to show that
the continued fraction of $\beta'$ is equal to
\begin{equation}\label{eq:102}
\beta'=[0,1,c,b_2,b_3,\dots]\,.
\end{equation}
The following lemma is crucial for construction of factors of a Sturmian word with maximal index.

\begin{lem}\label{l}
Let $u$ be a Sturmian word with slope $\beta$ having the continued
fraction $\beta=[0,1,b_2,b_3,\dots]$. Let $w\in\L(u)$, and let
$r\in\Q$, $r\geq 2$ be such that $v=w^r\in\L(u)$. Denote
$$
w'=\varphi(w)\qquad\hbox{and}\qquad v'=\varphi(v)A^c\,,
$$
where $\varphi$ is the morphism given by~\eqref{eq:101}. Then $v'$
is a rational power of $w'$ in a Sturmian word $u'$ with slope
$\beta'=[0,1,c,b_2,b_3,\dots]$.
\end{lem}

\pf If $|w|=1$, then necessarily $w=A$, $v=A^r$ for $2<r\leq
b_2+1$, $\varphi(w)=A^cB$, and $\varphi(v)A^c=(A^cB)^rA^c$ is a
factor of $u'$, since a Sturmian word with slope
$\beta'=[0,1,c,b_2,b_3,\dots]$ has blocks $A^c$, $A^{c+1}$
separated by single letters $B$.

If $|w|\geq 2$, let us write $w=w_1w_2$ so that $w_2\neq\epsilon$ and $v=(w_1w_2)^{\lfloor r\rfloor}w_1$.
Then $\varphi(v)A^c = \varphi(w^{\lfloor
r\rfloor})\varphi(w_1)A^c$. In order to show that $\varphi(v)A^c$
is a power of $\varphi(w)$, it suffices to show that
$\varphi(w_1)A^c$ is a prefix of $\varphi(w)$ or
$\varphi(w)\varphi(w)$. If $w_2$ starts with $A$ or $BA$, then
$\varphi(w_2)$ has prefix $A^c$ and thus $\varphi(w_1)A^c$ is a
prefix of $\varphi(w)=\varphi(w_1)\varphi(w_2)$. Since
$BB\notin\L(u)$, it remains to discuss the special case when
$w_2=B$. As $|w|\geq2$, we have $w_1\neq\epsilon$. Since
$w_2w_1\in\L(u)$, the word $w_1$ must start with the letter $A$
and therefore $\varphi(w_1)A^c$ is a prefix of
$\varphi(w_1B)\varphi(w_1B)=\bigl(\varphi(w)\bigr)^2$. \pfk

\begin{thm}
Let $u$ be a Sturmian word with slope
$\alpha=[0,1,a_2,a_3,\dots]$. Then for every $N\in\N$ there exists
a factor $w\in\L(u)$ with index at least equal to
$2+a_{N+1}+\frac{q_{N-1}-2}{q_N}$, where $q_N$ is the denominator
of the $N$-th convergent of $\alpha$.
\end{thm}

\pf For $N=1$ it follows from the continued fraction of $\alpha$
that $q_1=1$, $q_0=1$ and therefore we have to find a factor with
index $2+a_2-1=a_2+1$. It suffices to put $w=A$. Therefore we
consider $N\geq 2$. We shall construct the desired factor $w$ and
its power $v$ by $(N-1)$-fold application of Lemma~\ref{l}.
Consider the irrational number $\alpha_0$ with the continued
fraction $\alpha_0=[0,1,a_{N+1},a_{N+2,\dots}]$. Take a Sturmian
word $u^{(0)}$ with slope ${\alpha_0}$ and its factors
$w^{(0)}:=A$, $v^{(0)}:=A^{1+a_{N+1}}$ for initial values of the
construction. For $1\leq i\leq N-1$, define
$$
w^{(i)}:=\varphi_i(w^{(i-1)})\,,\qquad
v^{(i)}:=\varphi_i(v^{(i-1)})A^{a_{N-i+1}}\,,\qquad\hbox{where}\quad
\varphi_i:\begin{array}{rcl}A&\mapsto& A^{a_{N-i+1}}B\\ B&\mapsto&
A\end{array}\,.
$$
By Lemma~\ref{l}, the word $w^{(i)}$ is a factor of a Sturmian
word $u^{(i)}$ with slope ${\alpha_i}$, where $\alpha_i$ has the
continued fraction $\alpha_i=[0,1,a_{N+1-i}, a_{N+2-i},\dots]$ and
$v^{(i)}$ is a power of $w^{(i)}$ in the word $u^{(i)}$. In
particular, $w^{(N-1)}$ is a factor of a Sturmian word $u$ with
slope $\alpha=[0,1,a_2,a_3,\dots]$ and $v^{(N-1)}$ is its power in
$u$.

 It suffices now to show that the length of $w^{(N-1)}$ is $q_N$
and the length of $v^{(N-1)}$ is $(2+a_{N+1})q_N+q_{N-1}-2$. For the
recurrent expression of lengths of factors $w^{(i)}$, $v^{(i)}$ we use
formula~\eqref{eq:p3}. We have
$$
\bigl(|w^{(i)}|_A,|w^{(i)}|_B\bigr) =
\bigl(|w^{(i-1)}|_A,|w^{(i-1)}|_B\bigr) M_{a_{N-i+1}}\,,
$$
for all $i=1,2,\dots,N-1$, with
$\bigl(|w^{(0)}|_A,|w^{(0)}|_B\bigr)=(1,0)$. It can be easily seen
that
$$
\bigl(|w^{(N-1)}|_A,|w^{(N-1)}|_B\bigr) =
(1,0)M_{a_{N}}M_{a_{N-1}}\cdots M_{a_{2}}\,.
$$
In order to obtain $|w^{(N-1)}|=|w^{(N-1)}|_A+|w^{(N-1)}|_B$, we
multiply the latter from the right by the vector $\binom{1}{1}$,
which can be also written as
$\binom{1}{1}=\bigl(\begin{smallmatrix}1&1\\1&0\end{smallmatrix}\bigr)\binom{1}{0}$.
Since in the continued fraction of $\alpha$ we have $a_1=1$, we can
use~\eqref{eq:01} to obtain
$$
|w^{(N-1)}|=(1,0)M_{a_{N}}M_{a_{N-1}}\cdots
M_{a_{2}}M_{a_1}\textstyle{\binom{1}{0}}=q_N\,.
$$

From the definition of words $v^{(i)}$ we have for their lengths
\begin{equation}\label{eq:02}
\bigl(|v^{(i)}|_A,|v^{(i)}|_B\bigr) =
\bigl(|v^{(i-1)}|_A,|v^{(i-1)}|_B\bigr)
M_{a_{N-i+1}}+(a_{N-i+1},0)\,,
\end{equation}
with $\bigl(|v^{(0)}|_A,|v^{(0)}|_B\bigr)=(1+a_{N+1},0)$. Let us
compute the lengths for $N=1$,
$$
\bigl(|v^{(1)}|_A,|v^{(1)}|_B\bigr) =
(1+a_{N+1},0)\bigl(\begin{smallmatrix}a_N&1\\1&0\end{smallmatrix}\bigr)+(a_{N},0)=
(2+a_{N+1})(1,0)M_{a_{N}}+(1,0)-(1,1)\,.
$$
Since for every $c$ we have $-(1,1)M_c+(c,0)=-(1,1)$, by repeated
application of the recurrence~\eqref{eq:02} we obtain
$$
\bigl(|v^{(N-1)}|_A,|v^{(N-1)}|_B\bigr)=(2+a_{N+1})(1,0)M_{a_{N}}M_{a_{N-1}}\cdots
M_{a_{2}}+(1,0)M_{a_{N-1}}\cdots M_{a_{2}}-(1,1)\,.
$$
Again, multiplying the latter from the right by the vector
$\binom{1}{1}=M_{a_1}\binom{1}{0}$ and using~\eqref{eq:01}, we
obtain
$$
|v^{(N-1)}|=(2+a_{N+1})q_N+q_{N-1}-2\,.
$$
\pfk


\section{Acknowledgements}

We are grateful to J.-P. Allouche for pointing out the
reference~\cite{Carpi}. We also acknowledge financial support by
the grants MSM6840770039 and LC06002  of the Ministry of
Education, Youth, and Sports of the Czech Republic.



\end{document}